\documentclass[]{siamart171218}

\usepackage{algorithm}
\usepackage{algpseudocode}
\usepackage{amssymb}
\usepackage{bbm}
\usepackage{commath}
\usepackage{cleveref}
\usepackage{etoolbox}
\usepackage{tikz}
\usetikzlibrary{angles,calc,quotes}
\usepackage{pgfplots}
\pgfplotsset{compat=1.16}
\usepackage{xstring}

\newcommand{\matr}[1]{\ensuremath{\mathbf{#1}}}
\newcommand{\vectr}[1]{\ensuremath{\mathbf{#1}}}
\newcommand{\unitvector}[1]{\ensuremath{\boldsymbol{\hat{#1}}}}

\newcommand{\realset}{\ensuremath{\mathbb{R}}}
\newcommand{\bigo}[1]{\ensuremath{\mathcal{O}\del{#1}}}
\newcommand{\basis}[1]{\ensuremath{\vectr{\hat{e}_{#1}}}}
\newcommand{\dotproduct}[2]{#1 \cdot #2}
\newcommand{\stduniform}{\mathcal{U}}
\newcommand{\gammafun}[1]{\Gamma\del{#1}}
\newcommand{\betafun}[3][1]{\mathrm{B}\del{\ifthenelse{\equal{#1}{1}}{}{#1;}#2,#3}}
\newcommand{\normalizedbetafun}[3][1]{\mathrm{I}\del{\ifthenelse{\equal{#1}{1}}{}{#1;}#2,#3}}
\newcommand{\invnormalizedbetafun}[3][1]{\mathrm{I}^{-1}\del{\ifthenelse{\equal{#1}{1}}{}{#1;}#2,#3}}
\newcommand{\indicator}[1]{\ensuremath{\mathbbm{1}_{#1}}}
\newcommand{\anglemapsymbol}{\ensuremath{\Theta}}
\newcommand{\invanglemapsymbol}{\ensuremath{\anglemapsymbol^{-1}}}
\newcommand{\anglemap}[1]{\ensuremath{\anglemapsymbol\del{#1}}}
\newcommand{\invanglemap}[1]{\ensuremath{\invanglemapsymbol\del{#1}}}
\newcommand{\gaussianpdf}[1]{\ensuremath{\phi\del{#1}}}

\newrobustcmd{\ith}[1]{\StrLen{#1}[\stringlength]\ifnumequal{\stringlength}{1}{$#1$}{$(#1)$}-th}

\newcommand{\rejectioncostplanaranglea}{\ensuremath{\frac{\pi}{3}}}
\newcommand{\rejectioncostplanarangleb}{\ensuremath{\frac{\pi}{4}}}
\newcommand{\rejectioncostplanaranglec}{\ensuremath{\frac{\pi}{5}}}

\newcommand{\thetarejectioncostplanaranglea}{\ensuremath{\frac{\pi}{3}}}
\newcommand{\thetarejectioncostplanarangleb}{\ensuremath{\frac{\pi}{4}}}
\newcommand{\thetarejectioncostplanaranglec}{\ensuremath{\frac{\pi}{5}}}

\newcommand{\histogrambins}{100}
\newcommand{\histogramsamples}{10000}
\newcommand{\histogramdimension}{10}
\newcommand{\histogrammaxtheta}{\ensuremath{\frac{\pi}{4}}}

\newcommand{\conesamplingksdimensiona}{10}
\newcommand{\conesamplingksdimensionb}{100}
\newcommand{\conesamplingksdimensionc}{1000}

\newcommand{\shiftedsphereksmeanmag}{1.0}
\newcommand{\shiftedsphereksdimension}{10}

\newcommand{\reweightednormalksdimension}{100}
\newcommand{\reweightednormalksmeanmag}{1.0}
\newcommand{\reweightednormalksstandarddeviationa}{0.08}
\newcommand{\reweightednormalksstandarddeviationb}{0.12}
\newcommand{\reweightednormalacceptanceratioa}{0.9831}
\newcommand{\reweightednormalacceptanceratiob}{0.0968}

\def\power10typesetter#1{%
  \pgfkeys{/pgf/number format/.cd,sci,retain unit mantissa=false}%
  \pgfmathprintnumber{#1}%
}
\title{An $\mathcal{O}$(\lowercase{n}) algorithm for generating
  uniform random vectors in \lowercase{n}-dimensional cones}
\author{Arun I. \and Murugesan Venkatapathi
\thanks{Department of
    Computational and Data Sciences, Indian Institute of Science,
    Bengaluru - 560012. \email{murugesh@iisc.ac.in}}}

\begin{document}
\newcommand{\mudir}{\ensuremath{\unitvector{\mu}}}
\newcommand{\definestduniform}[1][$\stduniform$]{#1 is the standard
  uniform random variable distributed uniformly between 0 and 1}
\newcommand{\universalityoftheuniform}{Using the universality of the
  uniform random variable (also known as the probability integral
  transform), if \definestduniform, then}
\newcommand{\arccolor}{blue}
\newcommand{\arc}{{\color{\arccolor}arc}}
\newcommand{\sectorcap}[3][1]{\draw [\arccolor,very thick] ({(#2)-(#3)}:#1)
  arc [start angle=(#2)-(#3), end angle=(#2)+(#3), radius=#1]}
\newcommand{\sectorwithcap}[3][1]{\filldraw [green!20!white,
  draw=green!50!black] (0,0) -- ({(#2)-(#3)}:{(#1)}) arc [start angle=(#2)-(#3), end
  angle=(#2)+(#3), radius=#1] -- cycle;
  \sectorcap[#1]{#2}{#3};
}

\maketitle

\begin{abstract}
  Unbiased random vectors distributed uniformly in $n$-dimensional
  space are widely used, and the computational cost of generating a
  vector increases only linearly with $n$. On the other hand,
  generating uniformly distributed random vectors in its subspaces
  typically involves the inefficiency of rejecting vectors falling
  outside, or re-weighting a non-uniformly distributed set of
  samples. Both approaches become severely ineffective as $n$
  increases. We present an efficient algorithm to generate uniformly
  distributed random directions in $n$-dimensional cones, to aid
  searching and sampling tasks in high dimensions.
\end{abstract}

\section{Introduction}

The problem of generating unbiased random vectors appears widely, and
as described in \cref{sec:non-uniform-samples-in-a-sphere}, has a
relatively trivial solution that scales as \bigo{n} arithmetic
operations where $n$ is the dimension of the space
\cite{mode1992random}. This problem can be reduced to an accumulation
of random points with a uniform probability density on the surface of
the unit sphere. This description using the surface of the unit
sphere, allows us to effectively use geometry in describing the
algorithms proposed for uniformly sampling the region of interest
given by a part of the unit sphere. It can play a critical role in
searching, learning, and sampling tasks in high dimensions
\cite{basri2011neighbor, buchta2012k-means, parsons2004clustering}.

We may simply generate samples that are uniformly distributed on the
entire unit sphere, but only accept those that are within a
region. This gives us our desired uniform distribution within the
region of interest, and the number of rejections is determined by the
fraction of the surface (solid angle fraction) of the sphere that we
wish to sample. For a region on the unit sphere bound using the planar
angles between the position vectors, one observes that the fraction of
the total solid angle represented by this region rapidly decreases
with the dimension. This makes rejection sampling prohibitively
expensive in high dimensions and we demonstrate this quantitatively in
\cref{sec:rejection-sampling-cost}. We may also re-weight a
non-uniformly distributed set of samples if the probability
distribution is known. But, the re-weighting errors are known to
increase significantly in higher
dimensions\cite{diakonikolas2019robust,strohmer2000numerical} and we
demonstrate this with examples in \cref{sec:re-weighting}.

Note that many such naive methods are not effective for large $n$ as
they do not generate the required uniform distribution, or do so at a
prohibitive cost. While other preferred methods such as
Markov-Chain-Monte-Carlo (MCMC) can be significantly more efficient
than the above naive methods in generating uniformly distributed
points in an arbitrary volume, they may nevertheless scale as poorly
as \bigo{n^5} in the required computing effort
\cite{dyer1988complexity,
  kannan1997random,lovasz2003simulated}. \bigo{n^3} and \bigo{n^2}
methods that use linear transformations for uniformly sampling certain
regularly shaped surfaces are described in the literature
\cite{devroye2006nonuniform}. In this work, we present an \bigo{n}
method that uses a non-linear transformation to generate random points
uniformly distributed on a section of the surface of the sphere.

\section{Problem Statement}

\begin{figure} 
  \centering
  \newcommand{\coneposition}{30}
  \newcommand{\conewidth}{20}
  \newcommand{\sphericalcap}{{\color{\arccolor}spherical cap}
    \tikz[scale=0.5] \sectorcap{\coneposition}{\conewidth};}
  \newcommand{\cone}{cone \tikz[scale=0.5]
    \sectorwithcap{\coneposition}{\conewidth};}
  \begin{tikzpicture}[scale=2]
    \coordinate (O) at (0,0);
    \coordinate (conecenter) at (\coneposition:1.0);
    \coordinate (coneedge) at (\coneposition+\conewidth:1.0);
    \draw (O) circle (1);
    \sectorwithcap{\coneposition}{\conewidth};
    \draw [->] (O) -- (\coneposition:1.0) node [right] {\mudir};
    \draw pic ["$\theta$", draw=black,angle radius=1cm, angle
    eccentricity=0.8] {angle = conecenter--O--coneedge};
  \end{tikzpicture}
  \caption[Cone sampling problem statement]{It is required to generate
    points uniformly distributed on the \sphericalcap{}. The
    \sphericalcap{} is the set of all unit vectors that fall within
    some angle $\theta$ of the central axis \mudir.}
  \label{fig:problem-statement}
\end{figure}
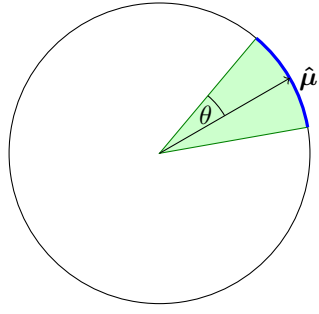

It is required to generate random points uniformly distributed on a
fraction $\Omega_0$ of the total solid angle $s_n$ of the unit
sphere. In two dimensions, this corresponds to generating random
points uniformly on an arc of the unit circle. In general, for $n$
dimensions, this corresponds to generating unit vectors in an
$n$-dimensional cone with a spherical cap at its base and its apex at
the center of the unit sphere of reference; we denote the central axis
of the cone from the apex to the centre of the base as \mudir{} (see
\cref{fig:problem-statement}). In other words, the desired cone is the
set of all unit vectors that fall within a planar angle $\theta$ of
the central axis \mudir{}. If $\mathrm{S}^{n-1}$ is the set of all
unit vectors, then this set is given by
\begin{equation}
  \set{\unitvector{x} \in \mathrm{S}^{n-1} ; \dotproduct{\unitvector{x}}{\mudir{}} \ge
    \cos \theta}
\end{equation}

where $0 \le \theta \le \pi$. Note that the proposed solution also extends to vectors contained in a
hollow cone bound by two planar angles.

\section{Map from planar angle to solid angle fraction}

If $\anglemapsymbol\colon\realset \to \realset$ is the function
mapping planar cross sectional angle $\theta$ to the rotated solid
angle fraction in $n$ dimensions,
$\phi, \theta_1, \theta_2, \dotsc, \theta_{n-2}$ are the $n - 1$
angles of the spherical coordinate system, and $s_n$ is the surface
area of the sphere, then \anglemap{\theta} is given by the following
integral where all angles except $\theta_{n-2}$ are integrated over
their full range, whereas $\theta_{n-2}$ is integrated over
\intcc{0,\theta}.
\begin{equation}
  \anglemap{\theta} = \frac{1}{s_n} \int_0^{2\pi} \int_0^{\theta} \int_0^{\pi} \dotsi
  \int_0^{\pi} \int_0^{\pi} \sin^{n-2}\theta_{n-2}
  \sin^{n-3}\theta_{n-3} \dotsm \sin^2\theta_2 \sin\theta_1 \dif \theta_{n-2} \dif \theta_{n-3} \dotsm \dif \theta_2 \dif \theta_1 \dif \phi
\end{equation}
The surface area $s_n$ of the sphere is given by
\begin{equation}
  s_n = \frac{2\pi^{\frac{n}{2}}}{\gammafun{\frac{n}{2}}}
\end{equation}
Separating out the multiple integral into a product of one-dimensional
integrals over each angle,
\newcommand{\angleintegral}[2][\pi]{\int_0^{#1} \sin{\ifthenelse{\equal{#2}{1}}{}{^{#2}}} \theta_{#2} \dif \theta_{#2}}
\begin{equation}
  \anglemap{\theta} = \frac{1}{s_n} \angleintegral[\theta]{n-2} \angleintegral{n-3} \dotsm
  \angleintegral{2} \angleintegral{1} \int_0^{2\pi} \dif \phi
\end{equation}
All but the first and last one of these integrals are Wallis'
integrals $W_m$ with $m = 1,2,3,\dotsc,n-3$.
\begin{equation}
  \int_0^{\pi} \sin^mx \dif x = 2 \int_0^{\frac{\pi}{2}} \sin^mx \dif
  x = 2 W_m
\end{equation}
Substituting $ u = \sin^2x $, we obtain the following relation between
the Wallis' integral and the beta function.
\begin{equation}
  W_m = \int_0^{\frac{\pi}{2}} \sin^mx \dif x = \frac{1}{2} \betafun{\frac{m+1}{2}}{\frac{1}{2}}
\end{equation}
Thus,
\begin{equation}
  \int_0^{\pi} \sin^mx \dif x = \betafun{\frac{m+1}{2}}{\frac{1}{2}}
\end{equation}
Substituting this integral back into \anglemap{\theta},
\newcommand{\betacomponent}[1]{B\del{\frac{#1}{2}, \frac{1}{2}}}
\begin{equation}
  \anglemap{\theta} = \frac{1}{s_n} 2\pi \betacomponent{n-2} \dotsm \betacomponent{3}
  \betacomponent{2} \angleintegral[\theta]{n-2}
\end{equation}
Expanding the beta function in terms of the gamma function
$\betafun{x}{y} = \frac{\gammafun{x}\gammafun{y}}{\gammafun{x+y}}$,
the product of beta functions telescopically cancels out to give
\begin{equation}
  \anglemap{\theta} = \frac{1}{\betafun{\frac{n-1}{2}}{\frac{1}{2}}} \angleintegral[\theta]{n-2}
\end{equation}
We can relate \anglemap{\theta} to the normalized incomplete Beta
function defined as
\begin{equation}
  \normalizedbetafun[x]{\alpha}{\beta} =
  \frac{1}{\betafun{\alpha}{\beta}} \int_0^x t^{\alpha-1}
  (1-t)^{\beta-1} \dif t
\end{equation}
In order to do so, we need to handle \anglemap{\theta} as two separate
cases---one when $\theta \in \intcc{0, \frac{\pi}{2}}$ and another
when $\theta \in \intcc{\frac{\pi}{2}, \pi}$.

\subsection{Case (i)}

When $\theta \in \intcc{0, \frac{\pi}{2}}$, the integral is an
incomplete form of the Wallis' integral, which we denote by
$W_m(\theta)$.
\begin{equation}
  W_m(\theta) = \int_0^{\theta} \sin^m x \dif x
\end{equation}
Substituting $u = \sin^2 x$, we obtain the following relation between
the incomplete Wallis' integral and the normalized incomplete Beta
function.
\begin{equation}
  W_m(\theta) = \frac{1}{2} \betafun{\frac{m+1}{2}}{\frac{1}{2}}
  \normalizedbetafun[\sin^2\theta]{\frac{m+1}{2}}{\frac{1}{2}} \quad \text{for} \enspace \theta \in \intcc{0, \frac{\pi}{2}}
\end{equation}

\subsection{Case (ii)}

When $\theta \in \intcc{\frac{\pi}{2}, \pi}$, we split the integral at
$\theta = \frac{\pi}{2}$.
\begin{equation}
  W_m(\theta) = \int_0^{\frac{\pi}{2}} \sin^m x \dif x +
  \int_{\frac{\pi}{2}}^{\theta} \sin^m x \dif x
\end{equation}
The first integral is the complete Wallis' integral $W_m$.
\begin{equation}
  W_m(\theta) = W_m + \int_{\frac{\pi}{2}}^{\theta} \sin^m x \dif x
\end{equation}
Due to the symmetry of the integrand $\sin^m x$ about
$\theta = \frac{\pi}{2}$,
\begin{equation}
  \int_{\frac{\pi}{2}}^{\theta} \sin^m x \dif x =
  \int_{\pi - \theta}^{\frac{\pi}{2}} \sin^m x \dif x =
  W_m - W_m(\pi - \theta)
\end{equation}
Thus,
\begin{equation}
  W_m(\theta) = 2W_m - W_m(\pi - \theta)
\end{equation}
\begin{equation}
  W_m(\theta) = \betafun{\frac{m+1}{2}}{\frac{1}{2}} \sbr{1 - \frac{1}{2}
    \normalizedbetafun[\sin^2\theta]{\frac{m+1}{2}}{\frac{1}{2}}} \quad \text{for} \enspace \theta \in \intcc{\frac{\pi}{2}, \pi}
\end{equation}
Thus,
\begin{equation}
  \int_0^{\theta} \sin^m x \dif x = \betafun{\frac{m+1}{2}}{\frac{1}{2}}
  \begin{cases}
    \frac{1}{2} \normalizedbetafun[\sin^2\theta]{\frac{m+1}{2}}{\frac{1}{2}} & \theta \in
    \intcc{0, \frac{\pi}{2}} \\
    1 - \frac{1}{2}
    \normalizedbetafun[\sin^2\theta]{\frac{m+1}{2}}{\frac{1}{2}} & \theta \in
    \intcc{\frac{\pi}{2}, \pi}
  \end{cases}
\end{equation}
And, the complete expression for \anglemap{\theta} is
\begin{equation}
  \anglemap{\theta} =
  \begin{cases}
    \frac{1}{2} \normalizedbetafun[\sin^2\theta]{\frac{n-1}{2}}{\frac{1}{2}} & \theta \in
    \intcc{0, \frac{\pi}{2}} \\
    1 - \frac{1}{2}
    \normalizedbetafun[\sin^2\theta]{\frac{n-1}{2}}{\frac{1}{2}} & \theta \in
    \intcc{\frac{\pi}{2}, \pi}
  \end{cases}
\end{equation}
If $\invnormalizedbetafun[y]{\alpha}{\beta}$ is the inverse of the
normalized beta function, computing the $x$ for which
\begin{equation}
  y = \normalizedbetafun[x]{\alpha}{\beta}
\end{equation}
then, the inverse of \anglemap{\theta} is given by
\begin{equation}
  \invanglemap{\Omega} =
  \begin{cases}
    \arcsin
    \sqrt{\invnormalizedbetafun[2\Omega]{\frac{n-1}{2}}{\frac{1}{2}}}
    & \Omega \in \intcc{0, \frac{1}{2}} \\
    \pi - \arcsin
    \sqrt{\invnormalizedbetafun[2\cbr{1-\Omega}]{\frac{n-1}{2}}{\frac{1}{2}}}
    & \Omega \in \intcc{\frac{1}{2}, 1}
  \end{cases}
\end{equation}

\section{Cost of high-dimensional rejection sampling}
\label{sec:rejection-sampling-cost}

The number of samples required to produce one accepted sample follows
a geometric distribution with the probability of acceptance given by
\begin{equation}
  p = \anglemap{\theta} = \frac{1}{\betafun{\frac{n-1}{2}}{\frac{1}{2}}} \int_0^{\theta}
  \sin^{n-2}\theta_{n-2}\dif\theta_{n-2}
\end{equation}
Thus, the average number of samples required to produce one accepted
sample is
\begin{equation}
  \label{eq:rejection-sampling-cost}
  \frac{1}{p} = \frac{1}{\anglemap{\theta}} =
  \frac{\betafun{\frac{n-1}{2}}{\frac{1}{2}}}{\int_0^{\theta} \sin^{n-2}\theta_{n-2}\dif\theta_{n-2}}
\end{equation}
This average number of samples, for various planar angles $\theta$, is
shown plotted against dimension $n$ in
\cref{fig:rejection-sampling-cost}. As can be seen from the figure,
the average number of samples required increases exponentially with
dimension. We justify this further using the following analytical
approximation. For a small planar angle $\theta$, $\theta_{n-2}$ only
takes on small values, and hence
$\sin\theta_{n-2} \approx \theta_{n-2}$. Therefore,
\begin{equation}
  \frac{1}{p} \approx
  \frac{\betafun{\frac{n-1}{2}}{\frac{1}{2}}}{\int_0^{\theta} \cbr{\theta_{n-2}}^{n-2} \dif\theta_{n-2}}
\end{equation}
Integrating and applying the limits,
\begin{equation}
  \frac{1}{p} \approx \frac{(n-1)\betafun{\frac{n-1}{2}}{\frac{1}{2}}}{\theta^{n-1}}
\end{equation}
For large $n$, we use Stirling's approximation for the beta function
\begin{equation}
  \betafun{\frac{n-1}{2}}{\frac{1}{2}} \approx \sqrt{\frac{2\pi
      e}{n-1}}
\end{equation}
to get
\begin{equation}
  \frac{1}{p} \approx \frac{\sqrt{2\pi e(n-1)}}{\theta^{n-1}}
\end{equation}
Thus, for small planar angles $\theta$, the average number of samples
required to produce one accepted sample increases exponentially with
dimension $n$. Therefore, rejection sampling is prohibitively
expensive especially in high dimensions. Note that while this
approximation applies only to small $\theta$, the cost of rejection
sampling rises exponentially even for large $\theta$. This follows
from \cref{eq:rejection-sampling-cost} and is observed in
\cref{fig:rejection-sampling-cost}.

\begin{figure}
  \centering
  \begin{tikzpicture}
    \begin{semilogyaxis}[
      enlarge x limits=false,
      enlarge y limits=false,
      cycle list={{black}},
      xlabel=Dimension, ylabel=Average number of samples]
      \addplot table {data/rejection-cost-pi-5.dat};
      \addplot table {data/rejection-cost-pi-4.dat};
      \addplot table {data/rejection-cost-pi-3.dat};
      \node at (axis cs:80,1e16) {$\theta=\rejectioncostplanaranglea$};
      \node at (axis cs:80,1e10) {$\theta=\rejectioncostplanarangleb$};
      \node at (axis cs:80,1e3) {$\theta=\rejectioncostplanaranglec$};
    \end{semilogyaxis}
  \end{tikzpicture}
  \caption{The average number of samples required by rejection
    sampling to produce one accepted sample, given a fixed planar
    angle, is shown plotted against dimension. As indicated in the
    plot, the three traces correspond to planar angles
    $\theta = \rejectioncostplanaranglea$,
    $\theta = \rejectioncostplanarangleb$ and
    $\theta = \rejectioncostplanaranglec$.}
  \label{fig:rejection-sampling-cost}
\end{figure}

\section{Proposed method of generation}

As a prerequisite, we will first look at generating points distributed
uniformly on a sphere. Then, we will describe the special case of
generating points uniformly distributed on the spherical cap of an
$n$-dimensional cone aligned along the \ith{n} canonical
axis. Finally, we will describe a rotation to reorient the random
points on the spherical cap to the desired direction \mudir{}.

\subsection{Generating random points on a sphere}
\label{sec:non-uniform-samples-in-a-sphere}

We can generate points uniformly distributed \emph{on} the surface of
the unit sphere using the Box-Muller transform \cite{box1958}. If
$Z_1, Z_2, \dotsc, Z_n$ are standard normal random variables and
$\basis{1}, \basis{2}, \dotsc, \basis{n}$ are canonical basis vectors,
then the vector \unitvector{s} uniformly distributed on the surface of
the sphere is given by
\begin{equation}
  \label{eq:box-muller}
  \unitvector{s} = \frac{Z_1\basis{1} + Z_2\basis{2} +
    Z_3\basis{3} + \dotsb + Z_n\basis{n}}{\sqrt{Z_1^2 +
      Z_2^2 + Z_3^2 + \dotsb + Z_n^2}}
\end{equation}

\subsection[Generating points along the nth canonical axis]{Generating
  points along the \ith{n} canonical axis}
\label{sec:generating-points-along-nth-canonical-axis}

\newcommand{\sphericalcap}{{\color{\arccolor}spherical cap}}
\newcommand{\diskcolor}{red}
\newcommand{\thetamax}{50}
\newcommand{\thetasample}{30}
\newcommand{\chord}{{\color{\diskcolor}chord}}
\newcommand{\disk}{{\color{\diskcolor}disk}}
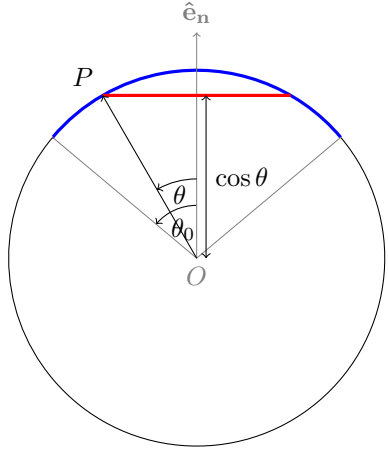
\newcommand{\geometrycaption}{It is required to generate random points
  P that are uniformly distributed on the \sphericalcap{}. The
  \sphericalcap{} subtends a fraction $\Omega_0 = \anglemap{\theta_0}$
  of the surface area of the sphere and is aligned along the \ith{n}
  canonical axis $\basis{n}$. To generate P, we generate samples
  uniformly distributed on the periphery of a \disk{} corresponding to
  some $\theta \in \intcc{0,\theta_0}$. $\theta$ is distributed such
  that samples on the disk uniformly sample the surface of the
  sphere.}
\begin{figure}
  \centering
  \begin{tikzpicture}[scale=2.5]
    \coordinate (O) at (0,0);
    \coordinate (C) at (90-\thetamax:1);
    \coordinate (D) at (90+\thetamax:1);
    \coordinate (P) at (90+\thetasample:1);
    \coordinate (Q) at (90-\thetasample:1);
    \coordinate (A) at (0,{cos(\thetasample)});
    % Draw circle, arc, nth canonical direction and mark sector angle \theta_0
    \draw [\arccolor,very thick] (C) arc [start angle=90-\thetamax, delta angle=2*\thetamax, radius=1] -- (D);
    \draw (C) arc [start angle=90-\thetamax, delta angle=-360+2*\thetamax, radius=1] -- (D);
    \draw pic [draw,->, angle radius=20, "$\theta_0$"] {angle = A--O--D};
    \draw [help lines] (D) -- (O) -- (C);
    \draw [->,help lines] (O) node [below] {$O$} -- (0,1.2) node [above] {$\basis{n}$};
    % Draw disk, and vector on it
    \draw [\diskcolor,very thick] (P) -- (Q);
    \draw [->] (O) -- (P) node [above left] {$P$};
    % Mark sin theta and cos theta
    \draw pic [draw,->, angle radius=30, angle eccentricity=0.83, "$\theta$"] {angle = A--O--P};
    % \draw [<->] ($(A) - (0,0.05)$) -- node [below] {$\sin \theta$} ($(P) - (0,0.05)$);
    \draw [<->] ($(O) + (0.05,0)$) -- node [right] {$\cos \theta$} ($(A) + (0.05,0)$);
  \end{tikzpicture}
  \caption{\geometrycaption}
  \label{fig:sampling-along-nth-canonical-axis}
\end{figure}

\geometrycaption{} Therefore, the cumulative distribution function of
$\theta$ is
\begin{equation}
  \label{eq:theta-cdf}
  F_{\theta}(\theta) =
  \begin{cases}
    \frac{\anglemap{\theta}}{\anglemap{\theta_0}} & 0 \le \theta \le \theta_0
    \\
    1 & \theta > \theta_0
  \end{cases}
\end{equation}
The probability density function of $\theta$ is
\begin{equation}
  \label{eq:theta-pdf}
  f_{\theta}(\theta) =
  \begin{cases}
    \frac{s_{n-1}}{s_n \anglemap{\theta_0}} \sin^{n-2}\theta & 0 \le
    \theta \le \theta_0 \\
    0 & \theta > \theta_0
  \end{cases}
\end{equation}
If $\theta$ is a random sample from this distribution,
$Z_1, Z_2, \dotsc, Z_{n-1}$ are standard normal random variables and
$\basis{1}, \basis{2}, \dotsc, \basis{n}$ are canonical basis vectors,
similar to \cref{eq:box-muller}, then we can construct a random vector
uniformly distributed on the spherical cap as
\begin{equation}
  \unitvector{x} = \sin \theta \frac{Z_1\basis{1} + Z_2\basis{2} +
    Z_3\basis{3} + \dotsb + Z_{n-1}\basis{n-1}}{\sqrt{Z_1^2 +
      Z_2^2 + Z_3^2 + \dotsb + Z_{n-1}^2}} + \cos \theta \basis{n}
\end{equation}
$\theta$ can be generated using inverse transform sampling as shown in
\cref{alg:generate-random-theta-inverse-transform-sampling}, or using
one-dimensional rejection sampling as shown in
\cref{alg:generate-random-theta-rejection-sampling}. Rejection
sampling to generate $\theta$ is, on the average, more expensive. But,
it avoids the need to compute \anglemapsymbol{} and
\invanglemapsymbol{}, both of which are vulnerable to floating point
underflow due to the small solid angle fractions involved when the
dimension $n$ is large and the angle $\theta_0$ is small. The $\log$
functions in \cref{alg:generate-random-theta-rejection-sampling} are
to further alleviate floating point underflow issues. Note that the
rejection sampling of
\cref{alg:generate-random-theta-rejection-sampling} is one-dimensional
and not subject to the prohibitive costs described in
\cref{sec:rejection-sampling-cost}. In fact, the cost of this
rejection sampling, that is, the average number of samples required to
produce a single $\theta$ is
\begin{equation}
  \frac{1}{\sqrt{\pi}}
  \frac{\Gamma\del{\frac{n}{2}}}{\Gamma\del{\frac{n-1}{2}}}
  \frac{\theta_0 \sin^{n-2}\cbr{\min(\theta_0,\frac{\pi}{2})}}{\Theta(\theta_0)}
\end{equation}
This cost, shown in \cref{fig:theta-rejection-sampling-cost} for
various $n$ and $\theta_0$, is approximately linear in $n$, and thus
the cost of the overall algorithm remains \bigo{n}.

\begin{figure}
  \centering
  \begin{tikzpicture}
    \begin{axis}[
      enlarge x limits=false,
      enlarge y limits=false,
      cycle list={{black}},
      xlabel=Dimension, ylabel=Average number of samples]
      \addplot table {data/theta-rejection-cost-pi-5.dat};
      \addplot table {data/theta-rejection-cost-pi-4.dat};
      \addplot table {data/theta-rejection-cost-pi-3.dat};
      \node at (axis cs:800,400) {$\theta_0=\thetarejectioncostplanaranglea$};
      \node at (axis cs:900,620) {$\theta_0=\thetarejectioncostplanarangleb$};
      \node at (axis cs:800,800) {$\theta_0=\thetarejectioncostplanaranglec$};
    \end{axis}
  \end{tikzpicture}
  \caption{The average number of samples required by rejection
    sampling to produce one accepted sample, given a fixed planar
    angle, is shown plotted against dimension. As indicated in the
    plot, the three traces correspond to planar angles
    $\theta = \thetarejectioncostplanaranglea$,
    $\theta = \thetarejectioncostplanarangleb$ and
    $\theta = \thetarejectioncostplanaranglec$.}
  \label{fig:theta-rejection-sampling-cost}
\end{figure}

This $\mathcal{O}(n)$ arithmetic procedure for generating random
vectors in a cone along the \ith{n} canonical axis \basis{n}, can be
extended to an arbitrary direction \mudir{} as shown in the next
section.

\subsection{Generating points along an arbitrary direction}
\label{sec:generating-points-along-arbitrary-direction}

We have generated points along the \ith{n} canonical axis \basis{n},
but it is required to generate points along a given arbitrary
direction \mudir{}. To do this, we simply rotate the vectors to align
along \mudir. More precisely, we rotate vectors $\hat{x}$ by the angle
between $\basis{n}$ and \mudir{} with the plane containing $\basis{n}$
and \mudir{} being the plane of rotation. If $\matr{P}$ is an
orthonormal matrix whose columns form a basis for the plane containing
$\basis{n}$ and \mudir, and $\matr{G}$ is the two dimensional Given's
rotation matrix for the required rotation, the rotated vector
$\hat{y}$ is
\begin{equation}
  \unitvector{y} = \unitvector{x} + \matr{P}\matr{G}\matr{P^T} \unitvector{x} - \matr{P}\matr{P^T}\unitvector{x}
\end{equation}
Rewriting using the $2 \times 2$ identity matrix $\matr{I_2}$,
\begin{equation}
  \unitvector{y} = \unitvector{x} + \matr{P}(\matr{G} - \matr{I_2})\matr{P^T}\unitvector{x}
\end{equation}
$\matr{P}$, $\matr{G}$ and $\mu_n$ are given by
\begin{equation}
  \matr{P} =
  \begin{bmatrix}
    \basis{n} & \frac{\mudir - \mu_n\basis{n}}{\norm{\mudir - \mu_n\basis{n}}}
  \end{bmatrix}
\end{equation}
\begin{equation}
  \matr{G} =
  \begin{bmatrix}
    \mu_n & - \sqrt{1-\mu_n^2} \\
    +\sqrt{1-\mu_n^2} & \mu_n
  \end{bmatrix}
\end{equation}
\begin{equation}
  \mu_n = \basis{n}^T \mudir
\end{equation}

A general $n$-dimensional rotation of vectors costs \bigo{n^2}
operations. But, the above is a simple\footnote{A simple rotation is a
  rotation with only one plane of rotation.} rotation and hence costs
only \bigo{n} operations. Thus, we have generated random points
aligned along an arbitrary direction without an increase in the order
of arithmetic complexity.

The overall algorithm discussed in
\cref{sec:generating-points-along-nth-canonical-axis,sec:non-uniform-samples-in-a-sphere,sec:generating-points-along-arbitrary-direction}
is shown in
\cref{alg:generate-random-point-on-sphere,alg:rotate-vector-from-nth-canonical-basis,alg:generate-random-theta-inverse-transform-sampling,alg:generate-random-theta-rejection-sampling,alg:generate-random-vector-on-spherical-cap}
where one of
\cref{alg:generate-random-theta-inverse-transform-sampling,alg:generate-random-theta-rejection-sampling}
can be used to generate random $\theta$.

\begin{algorithm}
  \begin{algorithmic}
    \Procedure{Generate point on sphere}{$n$}
    \For{$i = 1:n$}
    \State $x[i] \gets \text{normally distributed random number with
      (mean=0, variance=1)}$
    \EndFor
    \State $\vectr{x} \gets \frac{\vectr{x}}{\norm{\vectr{x}}}$
    \State \Return $\vectr{x}$
    \EndProcedure
  \end{algorithmic}
  \caption{Generate random point uniformly distributed on the
    $n$-dimensional unit sphere}
  \label{alg:generate-random-point-on-sphere}
\end{algorithm}
\begin{algorithm}
  \begin{algorithmic}
    \Procedure{Rotate vector from $n$-th canonical basis}{$\vectr{x}, \mudir$}
    \For{$i = 1:n-1$}
    \State $\matr{P}[i][1] \gets 0$
    \EndFor
    \State $\matr{P}[n][1] \gets 1 $
    \For{$i = 1:n-1$}
    \State $\matr{P}[i][2] \gets \frac{\hat{\mu}[i]}{\sqrt{1-\cbr{\hat{\mu}[n]}^2}}$
    \EndFor
    \State $\matr{P}[n][2] \gets 0 $
    \State $\matr{G} \gets \begin{bmatrix}
      \mu[n] & - \sqrt{1-\cbr{\mu[n]}^2} \\
      +\sqrt{1-\cbr{\mu[n]}^2} & \mu[n]
    \end{bmatrix}$
    \State $\matr{I_2} \gets \begin{bmatrix}
      1 & 0 \\
      0 & 1
    \end{bmatrix}$
    \State $\vectr{y} \gets \vectr{x} + \matr{P}(\matr{G} -
    \matr{I_2})\matr{P^T}\vectr{x}$
    \State \Return $\vectr{y}$
    \EndProcedure
  \end{algorithmic}
  \caption{Rotate vector \vectr{x} from around $n$-th canonical basis
    vector to arbitrary orientation \mudir{}}
  \label{alg:rotate-vector-from-nth-canonical-basis}
\end{algorithm}
\begin{algorithm}
  \begin{algorithmic}
    \Procedure{Generate random planar angle}{$\theta_0, n$}
    \State $\Omega_0 \gets \anglemapsymbol(\theta_0)$
    \State $U \gets \text{random number uniformly distributed between
      $0$ and $\Omega_0$}$
    \State $\theta \gets \invanglemap{U}$
    \State \Return $\theta$
    \EndProcedure
  \end{algorithmic}
  \caption{Generate random planar angle using inverse transform}
  \label{alg:generate-random-theta-inverse-transform-sampling}
\end{algorithm}
\begin{algorithm}
  \begin{algorithmic}
    \Procedure{Generate random planar angle}{$\theta_0, n$}
    \State $h \gets (n-2)\log \sbr{\sin \cbr{\min(\theta_0,
        \frac{\pi}{2})}}$
    \Repeat
    \State $U \gets \text{random number uniformly distributed between
      $0$ and $1$}$
    \State $\theta \gets \text{random number uniformly distributed between
      $0$ and $\theta_0$}$
    \State $f \gets h + \log U$
    \Until{$f < (n-2)\log\del{\sin \theta}$}
    \State \Return $\theta$
    \EndProcedure
  \end{algorithmic}
  \caption{Generate random planar angle using one-dimensional
    rejection sampling}
  \label{alg:generate-random-theta-rejection-sampling}
\end{algorithm}
\begin{algorithm}
  \begin{algorithmic}
    \Procedure{Generate point on spherical cap}{$\mudir,\theta_0$}
    \State $\theta \gets \text{Generate random planar angle}(\theta_0, n)$
    \State $\vectr{x} \gets \text{uninitialized $n$-dimensional vector}$
    \State $\vectr{x}[1:n-1] \gets \text{Generate point on sphere}(n-1)$
    \State $\vectr{x}[1:n-1] \gets \sin\theta \vectr{x}[1:n-1]$
    \State $\vectr{x}[n] \gets \cos \theta$
    \State $\vectr{x} \gets \text{Rotate vector from nth canonical
      basis}(\vectr{x}, \mudir)$
    \State \Return $\vectr{x}$
    \EndProcedure
  \end{algorithmic}
  \caption{Generate random points uniformly distributed on a spherical
    cap with axis \mudir{} and maximum planar angle $\theta_0$}
  \label{alg:generate-random-vector-on-spherical-cap}
\end{algorithm}

\subsection{Generating random vectors uniformly distributed in a hollow cone}

\newcommand{\hollowconeposition}{30}
\newcommand{\coneanglemin}{10}
\newcommand{\coneanglemax}{30}
\newcommand{\drawhollowcone}{
  \sectorwithcap{\hollowconeposition + \coneanglemin/2 +
    \coneanglemax/2}{\coneanglemax/2 - \coneanglemin/2};
  \sectorwithcap{\hollowconeposition - \coneanglemin/2 -
    \coneanglemax/2}{\coneanglemax/2 - \coneanglemin/2}
}
\newcommand{\hollowcone}{hollow cone \tikz[scale=0.5] {
    \drawhollowcone; }}

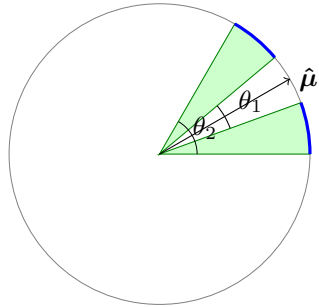
\begin{figure}[H]
  \centering
  \newcommand{\capofhollowcone}{{\color{blue}spherical cap} of the hollow cone \tikz[scale=0.5] { \drawhollowcone; }}
  \begin{tikzpicture}[scale=2]
    \coordinate (O) at (0,0);
    \draw [help lines] (O) circle (1);
    \drawhollowcone{};
    \coordinate (P) at (\hollowconeposition:1);
    \coordinate (A) at ({\hollowconeposition-\coneanglemax}:1);
    \coordinate (B) at ({\hollowconeposition+\coneanglemax}:1);
    \coordinate (C) at ({\hollowconeposition-\coneanglemin}:1);
    \coordinate (D) at ({\hollowconeposition+\coneanglemin}:1);
    \draw [very thin] pic ["$\theta_2$", draw=black, angle eccentricity=1.4] {angle = A--O--B};
    \draw pic ["$\theta_1$", draw=black,angle radius=1cm, angle eccentricity=1.4] {angle = C--O--D};
    \draw [->] (O) -- (\hollowconeposition:1) node [right] {\mudir};
  \end{tikzpicture}
  \caption[Hollow cone sampling problem statement]{It is required to
    generate points uniformly distributed on the \capofhollowcone{}
    along \mudir. The cross sectional planar angles $\theta_1$ and
    $\theta_2$ corresponding to the solid angles $\Omega_1$ and
    $\Omega_2$ are indicated.}
  \label{fig:hollow-cone-sampling-problem-statement}
\end{figure}

The proposed method can be generalized easily to generate points on
the surface of the unit sphere formed by the difference of two cones
corresponding to solid angle fractions $\Omega_1$ and $\Omega_2$ where
$\Omega_1 < \Omega_2$. Such a \hollowcone{} is shown in
\cref{fig:hollow-cone-sampling-problem-statement}. If
$\mathrm{S}^{n-1}$ is the set of all unit vectors, then the set of
unit vectors on the hollow cone is given by
\begin{equation}
  \set{\unitvector{x} \in \mathrm{S^{n-1}} ; \cos \theta_1 \ge
    \dotproduct{\unitvector{x}}{\mudir} \ge \cos \theta_2}
\end{equation}

If \definestduniform, the required random $\theta$ for such a
distribution of points is
\begin{equation}
  \theta = \invanglemap{\mathcal{U} (\Omega_2 - \Omega_1) + \Omega_1}
\end{equation}
$\theta_1$ and $\theta_2$ are the cross sectional planar angles
corresponding to solid angles fractions $\Omega_1$ and $\Omega_2$.
\begin{equation}
  \theta_1 = \invanglemap{\Omega_1}
\end{equation}
\begin{equation}
  \theta_2 = \invanglemap{\Omega_2}
\end{equation}

\section{Numerical validation}

\newcommand{\thetadefinition}{the angle between the generated
  direction vector and the reference direction \mudir{}}
\newcommand{\thetaodefinition}{The maximum $\theta$ is given by the
  cone angle $\theta_0=\histogrammaxtheta{}$}

The uniform distribution of the generated random points on the
spherical cap is verified by comparing the distribution of $\theta$,
\thetadefinition{}, with the exact analytically known distribution.\footnote{code freely available - https://pypi.org/project/sambal/ }

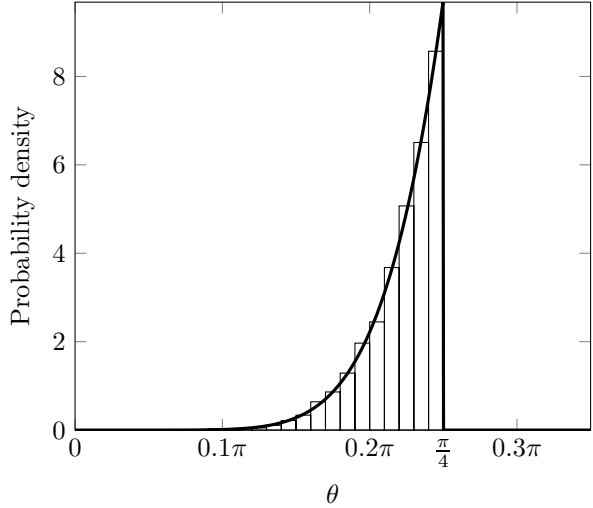
\newcommand{\histogramcaption}[1]{An empirical probability density of
  $\theta$ for a dimension $n=\histogramdimension{}$ and
  $\theta_0=\histogrammaxtheta{}$ is shown#1 with the exact
  probability density function overlaid. The empirical probability
  density was constructed using a histogram of \histogrambins{} bins
  and \histogramsamples{} direction vector samples.}
\histogramcaption{ in \cref{fig:histogram}} The exact probability
density function was given in \cref{eq:theta-pdf}.
\newcommand{\mytick}{\pgfmathparse{pi/4}\pgfresult}
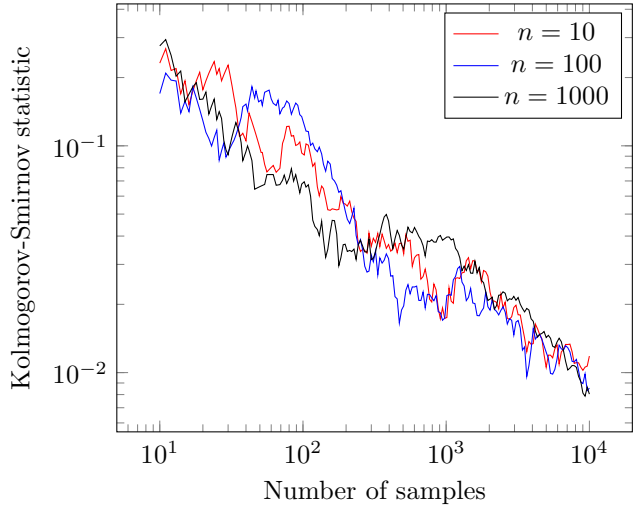
\newcommand{\kscaption}[1]{The Kolmogorov-Smirnov statistic comparing
  the empirical cumulative distribution function of $\theta$ and the
  exact cumulative distribution is plotted against an increasing
  number of samples#1. The Kolmogorov-Smirnov statistic decreases with
  an increasing number of samples indicating convergence of the
  empirical distribution to the exact distribution.}
\kscaption{ in \cref{fig:kolmogorov-smirnov}} The exact cumulative
distribution function was given in \cref{eq:theta-cdf}. If
$\theta_1,\theta_2,\dotsc,\theta_N$ are $N$ samples and
$\indicator{\intoc{-\infty,\theta}}$ is the indicator function of
$\intoc{-\infty,\theta}$, then the empirical cumulative distribution
function is
\begin{equation}
  F_N(\theta) = \frac{1}{N} \sum_{i=1}^N \indicator{\intoc{-\infty,\theta}}(\theta_i)
\end{equation}
The Kolmogorov-Smirnov statistic $D_N$ is the supremum of the absolute
difference between the empirical cumulative distribution function and
the exact cumulative distribution function.
\begin{equation}
  D_N = \sup_\theta \abs{F_N(\theta) - F(\theta)}
\end{equation}

\begin{figure}[H]
  \centering
  \begin{tikzpicture}
    \begin{axis}[
      enlarge x limits=false,
      enlarge y limits=false,
      xlabel=$\theta$,
      ylabel=Probability density,
      xtick={0.3*pi, 0.25*pi, 0.2*pi, 0.1*pi, 0},
      xticklabels={$0.3\pi$, $\frac{\pi}{4}$, $0.2\pi$, $0.1\pi$, $0$}]
      \addplot [ybar interval]
      table {data/histogram.dat};
      \addplot [black,very thick] table {data/histogram-exact.dat};
    \end{axis}
  \end{tikzpicture}
  \caption{\histogramcaption{} $\theta$ is \thetadefinition{}.}
  \label{fig:histogram}
\end{figure}

\begin{figure}[H]
  \centering
  \begin{tikzpicture}
    \begin{loglogaxis}[
      xlabel=Number of samples, ylabel=Kolmogorov-Smirnov statistic,
      cycle list name=color list,
      legend entries={$n=\conesamplingksdimensiona$, $n=\conesamplingksdimensionb$, $n=\conesamplingksdimensionc$}]
      \addplot table {data/cone-sampling-ks-dimension-\conesamplingksdimensiona.dat};
      \addplot table {data/cone-sampling-ks-dimension-\conesamplingksdimensionb.dat};
      \addplot table {data/cone-sampling-ks-dimension-\conesamplingksdimensionc.dat};
    \end{loglogaxis}
  \end{tikzpicture}
  \caption{\kscaption{} $\theta$ is
    \thetadefinition{}. \thetaodefinition{}.}
  \label{fig:kolmogorov-smirnov}
\end{figure}

\appendix

\section{A note on generating samples by re-weighting non-uniform
  distributions} \label{sec:re-weighting}

One could sample the desired region of interest using a non-uniform
distribution, where the samples are re-weighted for uniformity. Below
we present a few straight-forward methods to generate the samples
efficiently in a region of interest using \bigo{n} arithmetic
operations. When the non-uniform distribution of these samples is
known, they can be re-weighted to achieve a semblance of uniform
distribution in the region. However this approach produces poor
results compared to the proposed method, especially when $n$
increases.

\subsection{Shifting generated points to the regions of interest and
  re-weighting}

\begin{figure} [H]
  \centering
  \begin{tikzpicture}[scale=3]
    \coordinate (O) at (0,0);
    \coordinate (mu) at (30:1);
    \node [circle,draw,fill=green!20] (sphere) at (mu) [minimum
    size=50] {};
    \draw (O) circle (1);
    \coordinate (tangent) at (tangent cs:node=sphere,point={(O)},solution=1);
    \draw (O) -- (tangent) -- node [right] {$r_0$} (mu) -- node [above] {\norm{\vectr{\mu}}} cycle;
    \draw pic [draw,->, fill=green!20, angle radius=40, angle
    eccentricity=0.9, "$\theta_0$"] {angle = tangent--O--mu};
    \draw pic [draw,->, angle radius=5] {right angle = O--tangent--mu};
  \end{tikzpicture}
  \caption{We generate random points uniformly distribute in the
    smaller sphere and normalize the position vectors to get direction
    vectors that lie on the surface of the larger sphere. The smaller
    sphere is shifted from the origin by \mudir{}. The generated
    direction vectors lie within the cone of half-angle $\theta_0$.}
  \label{fig:shifted-sphere}
\end{figure}
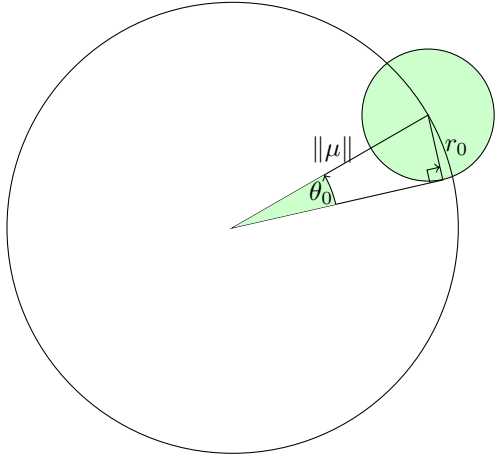

It is inexpensive to generate random points uniformly distributed in a
sphere. So, we bound our desired region of the unit sphere in a second
sphere, generate points inside that sphere and normalize the position
vectors of the generated points to get direction vectors that lie on
the surface of the first sphere. This geometry is shown graphically in
\cref{fig:shifted-sphere}. If $\vectr{S}$ is a vector random variable
distributed uniformly in the unit sphere centered at the origin, then
the generated direction vector \unitvector{x} is
\begin{equation}
  \vectr{x} = r_0\vectr{S} + \vectr{\mu}
\end{equation}
The radius $r_0$ depends on $\theta_0$ as follows.
\begin{equation}
  r_0 = \norm{\vectr{\mu}} \sin \theta_0
\end{equation}
\begin{equation}
  \unitvector{x} = \frac{\vectr{x}}{\norm{\vectr{x}}}
\end{equation}
The probability density of the generated direction vector at
\unitvector{x} is
\begin{equation}
  f(\unitvector{x}) = A \del{r_1^n - r_2^n}
\end{equation}
where $A$ is the normalization constant. $r_1$ and $r_2$ are
\begin{equation}
  r_1,r_2 = \del{\dotproduct{\unitvector{x}}{\vectr{\mu}}} \pm
  \sqrt{\del{\dotproduct{\unitvector{x}}{\vectr{\mu}}}^2 - \norm{\vectr{\mu}}^2 + r_0^2}
\end{equation}
Note that the choice of bounding the desired region in a sphere is
arbitrary. We could also have bound the desired region in a cube or
any other convenient shape.

A comparison of this method with the proposed method using the
Kolmogorov-Smirnov statistic of $\theta$ is shown in
\cref{fig:kolmogorov-smirnov-shifted-sphere}. $\theta$ is
\thetadefinition. If
$\unitvector{x}_1,\unitvector{x}_2,\dotsc,\unitvector{x}_N$ are $N$
samples and $\indicator{\intoc{-\infty,\theta}}$ is the indicator
function of $\intoc{-\infty,\theta}$, then the empirical cumulative
distribution function of $\theta$ after re-weighting is given by
\begin{equation}
  F_N(\theta) = \frac{1}{\sum_{i=1}^N \frac{1}{f(\unitvector{x}_i)}}
  \sum_{i=1}^N \frac{\indicator{\intoc{-\infty,\theta}}\del{\cos^{-1}\cbr{\dotproduct{\unitvector{x}}{\mudir}}}}{f(\unitvector{x}_i)}
\end{equation}

\begin{figure}[H]
  \centering
  \begin{tikzpicture}
    \begin{loglogaxis}[
      enlarge x limits=false,
      xlabel=Number of samples, ylabel=Kolmogorov-Smirnov statistic,
      cycle list name=color list]
      \addplot table
      {data/shifted-sphere-ks.dat}
      [yshift=24pt] node [above,pos=0.5] {Shifted sphere};
      \addplot table {data/cone-sampling-ks-dimension-\shiftedsphereksdimension.dat}
      [yshift=-24pt] node [below,pos=0.5] {Proposed method};
    \end{loglogaxis}
  \end{tikzpicture}
  \caption{The Kolmogorov-Smirnov statistic of the re-weighted
    distribution of shifted sphere random vectors for $\theta$ is
    compared with that of the proposed method for a dimension
    $n=\shiftedsphereksdimension$. $\theta$ is
    \thetadefinition{}. \thetaodefinition{}, and
    $\norm{\vectr{\mu}} = \shiftedsphereksmeanmag$. The larger value of $n$=100 is not shown here, as the re-weighted samples do not exhibit convergence.}
  \label{fig:kolmogorov-smirnov-shifted-sphere}
\end{figure}

\subsection{Re-weighting a shifted multivariate normal distribution}

If the probability density function of a non-uniformly distributed
random vector generator is known, the non-uniformly distributed
samples can be re-weighted to obtain a uniformly distributed random
vector in the limit of a large number of samples. However, error
introduced due to re-weighting can hinder performance. We demonstrate
this here with a multivariate normal distribution.

Generate vectors \vectr{x} that are distributed as the multivariate
normal distribution $\mathcal{N}(\vectr{\mu}, \sigma^2\matr{I_n})$
with mean \vectr{\mu} and covariance $\sigma^2\matr{I_n}$, where
\matr{I_n} is the $n \times n$ identity matrix and \vectr{\mu} is a
vector of arbitrary magnitude in the desired direction
\mudir{}. Normalize \vectr{x} to get the direction vector
\unitvector{x}. The probability density at \unitvector{x} is given by
\begin{equation}
  \label{eq:gaussian-pdf}
  f(\hat{x}) = \frac{\gaussianpdf{\frac{1}{\sigma}
      \sqrt{\norm{\vectr{\mu}}^2-(\dotproduct{\hat{x}}{\vectr{\mu}})^2}}}{\sigma^n
    (2\pi)^{\frac{n}{2} - 1}} \int_0^{\infty} r^{n-1}
  \gaussianpdf{\frac{r - \dotproduct{\hat{x}}{\vectr{\mu}}}{\sigma}} \dif r
\end{equation}
where $\phi$ is the probability density of the standard univariate
normal distribution.
\begin{equation}
  \gaussianpdf{x} = \frac{1}{\sqrt{2\pi}} e^{-\frac{x^2}{2}}
\end{equation}

If \unitvector{x} does not fall within the desired spherical cap,
reject and repeat. The probability density at \unitvector{x} after
rejection is proportional to the probability density of
\cref{eq:gaussian-pdf}. The value of $\sigma$ should be chosen to
minimize the rejection, but that makes the distribution more
non-uniform resulting in larger re-weighting errors. A comparison of
this method with the proposed method using the Kolmogorov-Smirnov
statistic is shown in \cref{fig:kolmogorov-smirnov-reweighted-normal}.

\begin{figure}[H]
  \centering
  \begin{tikzpicture}
    \begin{loglogaxis}[
      xmin=10, enlarge x limits=false,
      xlabel=Number of samples, ylabel=Kolmogorov-Smirnov statistic,
      cycle list name=color list]
      \addplot table
      {data/reweighted-normal-ks-standard-deviation-\reweightednormalksstandarddeviationa.dat}
      [yshift=28pt] node [above,pos=0.65] {Normal ($\sigma=\reweightednormalksstandarddeviationa$)};
      \addplot table
      {data/reweighted-normal-ks-standard-deviation-\reweightednormalksstandarddeviationb.dat}
      [yshift=0pt] node [above,pos=0.1] {Normal ($\sigma=\reweightednormalksstandarddeviationb$)};
      \addplot table {data/cone-sampling-ks-dimension-100.dat}
      [yshift=-20pt] node [below,pos=0.5] {Proposed method};
    \end{loglogaxis}
  \end{tikzpicture}
  \caption{The Kolmogorov-Smirnov statistic of a re-weighted
    multivariate normal distribution for $\theta$ is compared with
    that of the proposed method for a dimension
    $n=\reweightednormalksdimension$. $\theta$ is
    \thetadefinition{}. \thetaodefinition{}, and
    $\norm{\vectr{\mu}} = \reweightednormalksmeanmag$. For the
    re-weighted normal distributions, the fraction of accepted samples
    was \reweightednormalacceptanceratioa{} and
    \reweightednormalacceptanceratiob{} for
    $\sigma=\reweightednormalksstandarddeviationa$ and
    $\sigma=\reweightednormalksstandarddeviationb$ respectively.}
  \label{fig:kolmogorov-smirnov-reweighted-normal}
\end{figure}

\bibliographystyle{siamplain}
\bibliography{references}

\end{document}